\documentclass[11pt,reqno]{amsart}

 \usepackage{mathrsfs}
\newcommand{\scrL }{\mathscr{L}}

\newcommand{\scrS }{\mathscr{S}}

\usepackage{color}
\usepackage{amsmath, amsthm, amssymb}
\usepackage{amsfonts}
\usepackage[ansinew]{inputenc}
\usepackage[dvips]{epsfig}
\usepackage{graphicx}
\usepackage[english]{babel}
\usepackage{hyperref}

\usepackage{cite}
\usepackage{graphicx}
\usepackage{amscd}
\usepackage{color}
\usepackage{bm}
\usepackage{enumerate}

\usepackage{verbatim}
\usepackage{hyperref}
\usepackage{amstext}
\usepackage{latexsym}
\theoremstyle{plain}
\newtheorem{Theorem}{Theorem}[section]
\newtheorem{Definition}[Theorem]{Definition}

\newtheorem{Proposition}[Theorem]{Proposition}
\newtheorem{Lemma}[Theorem]{Lemma}
\newtheorem{Remark}[Theorem]{Remark}

\numberwithin{Theorem}{section}
\numberwithin{equation}{section}

\def\proof{\noindent{{\bf Proof. }}}
\def\square{\vbox{
\hrule height .4pt \hbox{\vrule width .4pt height 7pt \kern 7pt
\vrule width .4pt} \hrule height .4pt }}
\def\QED{\hfill {$\square$}\goodbreak \medskip}

\newcommand{\average}{{\mathchoice {\kern1ex\vcenter{\hrule height.4pt
width 6pt depth0pt} \kern-9.7pt} {\kern1ex\vcenter{\hrule
height.4pt width 4.3pt depth0pt} \kern-7pt} {} {} }}

\def\R{\mathbb{R}}

\def\div{\text{div}}

\renewcommand{\a }{\alpha }
\renewcommand{\b }{\beta }

\newcommand{\e }{\varepsilon }

\renewcommand{\l }{\lambda }

\newcommand{\n }{\nabla }
\newcommand{\vp }{\varphi }
\renewcommand{\phi}{\varphi}

\newcommand{\s }{\sigma }
\newcommand{\Sig }{\Sigma}
\renewcommand{\t }{\tau }
\newcommand{\z }{\zeta}
\renewcommand{\th }{\theta }

\renewcommand{\o }{\omega }
\renewcommand{\O }{\Omega }

\newcommand{\ov}{\overline}

\newcommand{\be}{\begin{equation}}
\newcommand{\ee}{\end{equation}}

\newcommand{\de}{\partial}

\newcommand{\ti}{\widetilde}

\newcommand{\ra}{{\rangle}}
\newcommand{\la}{{\langle}}

\renewcommand{\k}{\kappa}

\newcommand{\calF}{{\mathcal F}}

\newcommand{\N}{\mathbb{N}}

\newcommand{\Z}{\mathbb{Z}}

\newcommand{\cK}{{\mathcal K}}

\newcommand{\cM}{{\mathcal M}}

\newcommand{\cO}{{\mathcal O}}

\newcommand{\cR}{{\mathcal R}}

\newcommand{\cV}{{\mathcal V}}

\newcommand{\M}{\cM}

\renewcommand{\phi}{\varphi}

\renewcommand{\epsilon}{\varepsilon}

\renewcommand{\u}{{u}}    

\newcommand{\tauE}{\tau_{\mbox{\tiny $E$}}}

\begin{document}
\author[Mouhamed Moustapha Fall]
{Mouhamed Moustapha Fall}
\address{M. M. Fall: African Institute for Mathematical Sciences of Senegal, 
KM 2, Route de Joal, B.P. 14 18. Mbour, S\'en\'egal}
\email{mouhamed.m.fall@aims-senegal.org}
\title{Constant  Nonlocal Mean Curvatures surfaces and related problems }
\date{}
\thanks{The   author would like to thank   Tobias Weth, Xavier Cabr\'e and Enrico Valdinoci for useful discussions.  The  author's work is supported by the Alexander  von Humboldt foundation. }
\maketitle

\section{Introduction}
The concept of \textit{mean curvature}  of a surface goes back to Sophie Germain's work on elasticity theory in the    seventeenth century.  The mathematical formulation of the mean curvature was first derived by Young and then by Laplace in the eighteenth, see \cite{Finn}.   The mean curvature   of a surface   is an extrinsic measure of curvature which   locally describes the curvature of  surface in some ambient space.  The notion of \textit{nonlocal mean curvature}    appeared recently and for the  first time    in the work of  Caffarelli and Souganidis  \cite{Caff-Soug2010}   on   cellular automata. It is also an extrinsic geometric quantity that is invariant under global reparameterization of a surface. While Constant Mean Curvature (CMC) surfaces are stationary points for the area functional under some   constraints, Constant  Nonlocal Mean Curvature (CNMC) surfaces are also critical points of a fractional order area functional, called fractional or nonlocal perimeter, under some   constraints.  For simplicity, the fractional perimeter of a bounded set is given by a Sobolev fractional seminorm of its indicator function. Moreover, up to a normalization, as the fractional parameter tends to a maximal value, it approaches  the classical perimeter functional. 
In  some  phase transition problems driven by L\`evy-type diffusion, the sharp-interface energy is given by the fractional perimeter functional, see \cite{SO}. Moreover,    the rich  structure of CNMC surfaces captures the geometry and distributions of some periodic equilibrium patterns observed in  physical systems involving short  and long range competitions, see \cite{FALL-polymer}.

In this note, we survey recent results on embedded  surfaces with nonzero constant nonlocal mean curvature together with  their counterparts  within the classical   theory of constant mean curvature  surfaces such as the Alexandrov's classification theorem,  rotationally symmetric and periodic surfaces.  Unlike the study of CMC surfaces,  where  we can rely on   the theory   ordinary differential equations in certain cases,  the nonlocal case does not provide such advantage.  Neither  does it   provide, in general, an analytic verification of parameterized surfaces to have  CNMC.  Up to now, 	almost all nontrivial surfaces with  CNMC   appearing in the literature are  either derived by variational methods or by means of topological bifurcation theory.  The  study of solutions to partial differential equations with overdetermined   boundary values      is surprisingly intimately related with the question of finding CNMC surfaces. We also describe a number of such phenomena, notably those that were found recently as formal consequences of the highly nonlinear and nonlocal  equations which are involved. 
\section{Constant mean curvature hypersurfaces}
Let $\Sig$ be an orientable $C^2$ hypersurface of $\R^N$ and denote by  ${\cV}_{\Sig}: \Sig\to \R^N$   the unit   normal vector field on $\Sig$. 
 For every $p\in \Sig$, we let  $\{e_1;\dots;e_{N-1}\}$ be  an orthonormal basis of the tangent plane $T_p\Sig$ of $\Sig$ at $p$. 
The (normalized) \textit{mean curvature}   at $p$ of $\Sig$ is given by 
 \be \label{eq:def-MC-geom}
 H(\Sig;p):=\frac{1}{N-1}\sum_{i=1}^{N-1}  \la D \cV_\Sig (p)e_i, e_i \ra. 
\ee 
  Here and in the following, $\la ,\ra$ and "$\cdot$" denote  scalar product on $\R^N$.    As a consequence,  for a $C^1$-extension of ${\cV}_{\Sig} $  by a unit vector field  $\ti{\cV}_{\Sig}$     in a neighborhood of $p$ in $\R^N$, we have 
$$
 H(\Sig;p):=\frac{1}{N-1}\div_{\R^N} \ti{\cV}_{\Sig} (p). 
$$
Let $\O$ and $E$ be two open subsets of $\R^N$ with $E\subset \O$. Then the 
     \textit{perimeter functional of $E$ relative to $\O$}  (total variation of  $1_E$ in $\O$) is given by 
$$ 
P(E,\O)=|D 1_E|(\O):= \sup \left\{ \int_{ E}{\rm div} \xi(x)  \, d x\,:\,{\xi\in C^\infty_c(\O;\R^N)}, \, {|\xi |\leq 1}   \right\}.
$$
  In the following, we will simply write $P(E):=P(E,\R^N)$.
\begin{Definition}
We consider a vector field  $\zeta\in C^\infty_c(\R^N,\R^N)$ and define the flow   $(Y_t)_{t\in \R}$   induced by $\zeta$ given by 
\begin{align}
\begin{cases}
\de_t Y_t (x)= \zeta( Y_t(x))& \qquad \textrm{ $t\in\R$}\\
Y_0(x)= x& \qquad \textrm{ for all $x\in  \R^N$}.
\end{cases} 
\end{align}
For $E\subset \R^N $, we  call  the family of sets $E_t:=Y_t(E)$, $t\in \R$,  the variation of $E$ with respect to the vector field $\zeta$.  
\end{Definition}
Physically, constant   mean curvatures  surfaces    appear also when looking at surface at equilibrium enclosing a given volume in the absence of gravity.   We have the well know formula for first variation of area, see \cite{Spivak}.
\begin{Proposition} 
Let $\O$ and $E$ be two bounded    domains of $\R^N$, with  $E$ of class $C^2$. Let  $\l\in \R$ and $(E_t)_t$ be a variation of $E$ with respect to $\zeta\in C^\infty_c(\O;\R^N)$.   Then the map $t\mapsto J(t):= P(E_t,\O)-\l |E_t\cap\O|$ is differentiable at zero. 
Moreover 
$$
J'(0) =(N-1)\int_{\de E  } \left\{H(\de E; p)- \l  \right\} v(p)\, dV(p) ,
$$
where  $v(p):= \la \zeta(p),  {\cV}_{\de E}(p)   \ra$ and     ${\cV}_{\de E} $ is the unit exterior normal vector field  of $ E$.
\end{Proposition}

The first classification result of compact CMC surfaces is due to Jellet in the $18^{\textrm{th}}$ in \cite{Jellet}, showing that  a compact $CMC$  surface,  enclosing a  star-shaped domain
  in $\R^3$ must be  the standard sphere.   An other classification is due to Alexandrov in \cite{Alexandrov}, which we record in the following
\begin{Theorem} \label{th:Alexandrov}
An  embedded closed   $C^2$  hypersurface in $\R^N$, with nonzero constant mean curvature    is a finite union of disjoint round spheres with same radius.
\end{Theorem}
Alexandrov's result for embedded CMC hypersurfaces provides, in particular, a positive partial answer to a conjecture of H. Hopf, stating that: a compact orientable CMC surface must be a sphere.
Hopf gave a positive answer to his conjecture in \cite{Hopf} for the case of CMC immersions of $S^2$ into $\R^3$. As what concerns immersed hypersurfaces  with genus larger than 1,     Hsiang 
found in \cite{Hsiang} a first counterexample in $\R^4$ to Hopf's conjecture. Later Wente constructed in \cite{Wente}
an immersion of the torus $\mathbb{T}^2$ in $\R^3$ having constant mean curvature. While 
Wente's  CMC torus has genus $g =1$,   compact CMC immersions     with any genus  $g>1$ in $\R^3$ were
obtained   by N. Kapouleas \cite{Kapouleas}.  For a recent survey on constant mean curvature surfaces, we refer the reader to \cite{MPT} and the references therein.\\
Alexandrov introduced the method of moving plane in the proof of  Theorem \ref{th:Alexandrov}.  Formally, the argument  goes as follows. Consider a connected hypersurface $\Sig$ of class $C^2$ and pick a unit vector $e\in \R^N$ together with a hyperplane $P_e$ which is perpendicular to $e$ and not intersecting $\Sig$. 
Then slide the plane along   the $e$-direction, toward $\Sig$, until one of the two properties occurs for the first time: 
\begin{enumerate}
\item[1.] \textit{interior touching:} the reflection of $\Sig$ with respect to $P_e$, called $\Sig_e$, intersects $\Sig$ at a point $p_0\not \in P_e$,
\item[2.] \textit{non-transversal intersection:} the vector $e\in T_{p_0} \Sig$, for some $p_0\in \Sig\cap\Sig_e$.
\end{enumerate}
In either case, local  comparison   principles for elliptic equations applied to the  graphs  $u$ and $u_e$, locally representing $\Sig$ and $\Sig_e$, show that there exists an $r>0$ such that $\Sig\cap B_r(p_0)=\Sig_e\cap B_r(p_0)$. Now by unique continuations property of elliptic PDEs,  we find that  $\Sig=\Sig_e$. Since $e$ is arbitrary, this implies that $\Sig$ must be a round sphere. Indeed, letting $H$ the mean curvature of $\Sig$, then    $u$ and $u_e$ satisfies 
$$
\div\frac{\n  u}{\sqrt{1+|\n u|^2}}=\div \frac{\n u_e}{\sqrt{1+|\n u_e|^2}}=H \qquad\textrm{ on $B$},
$$
for some centered ball $B\subset \R^{N-1} \simeq T_{p_0} \Sig$. Moreover, making the ball smaller if necessary, we have $u\geq u_e$  on $B$.  The main point here is that  $w=u-u_e$  solves an elliptic equation of the form
$$
\de_i(a_{ij}(x)\de _j w)=0 \qquad\textrm{ on $B $},
$$
for a symmetric matrix $(a_{ij}(x))_{1\leq i,j\leq N-1}$ with positive and $C^1$ coefficients, only depending on the partial derivatives of $u$ and $u_e$. The comparison principles that is need is resumed in the following result.
\begin{Lemma}\label{lem:mmp}
Let $w\in C^2(B)\cap C^1(\ov{B })$ be a nonegative function on $\ov B$ and  satisfy
$$
 \de_i(a_{ij}(x)\de_j w)=0 \qquad\textrm{ on $B $},
$$
for some positive definite matrix $a_{ij}$ of class $C^1$.   Then the following   holds.
\begin{enumerate}
\item If $w(x_0)=0$, for some $x_0\in \de B$ then $w=0$ in $ B$.
\item If $w(x_0)=\n w(x_0)\cdot \nu=0$, for some $x_0\in \de B$ and  $\nu $ a unit vector normal to  $T_{x_0}\de B$,   then $w=0$ in $ B$.
\end{enumerate}
\end{Lemma}
The proof of   Lemma \ref{lem:mmp} can be found in many text books, see e.g. \cite{GT}. 
 Lemma \ref{lem:mmp}$(i)$ and $(ii)$ are  now used, to deal with cases \textit{1.} and \textit{2.}, respectively, yielding $w=0$ in $ B$.    
%

%

\subsection{Rotationally symmetric  constant mean curvature  surfaces}\label{s:Del-local}
An important class of CMC surfaces can be found in the family of unbounded rotationnaly symmetric ones. They were first studied in 1841 by Delaunay   in \cite{Delaunay}, which  he described explicitly as surfaces of revolution of \textit{roulettes of the conics}\footnote{   http://www.mathcurve.com/courbes2d/delaunay/delaunay.shtml}.
These surfaces are the \textit{catenoids, unduloids, nodoids and right circular cylinders}. They have the form
$$
 \Sig=\{ (x(s), y(s)\s)\,:\, s\in \R,\, \s\in S^1 \}\subset \R \times \R^2. 
$$
 where the rotating  plane curve  $s\mapsto  (x(s),y(s))$  is parmeterized by  arclength and $y$ a positive function on $\R$.   
 The explicit form of the mean curvature of $\Sigma$ at a point $ q(s,\s): = (x(s), y(s)\s)$ is given by
 \be \label{eq:NMC-surf-rev}
 H(\Sig;q(s,\s))=\frac{1}{2}\frac{-x'(s)+ y(s)\{x'(s)y''(s)-x''(s)y'(s)\}}{y(s)} .
 \ee
Of particular interest to us in this note are  the \textit{embedded} surfaces with nonzero CMC, the  so called \textit{unduloids}.  They constitute a smooth  1-parameter family of surfaces $(\Sigma_b)_{b\in (0,1)}$ varying from the straight cylinder $ \Sig_0$ to a translation invariant tangent spheres $\Sig_1$. Their explicit form was found by Kenmotsu in \cite{Kenmotsu},  %
$$
 \Sig_b=\{ (x_b(s), y_b(s)\s)\,:\, s\in \R,\, \s\in S^1 \}, 
$$
 where 
$$ 
x_b(s)= \int_{0}^s\frac{  1+b\sin(h r) }{\sqrt{1+b^2+2b \sin(hr)}}\, dr  \qquad \textrm{ and }  \qquad y_b(s)=\frac{1}{h}\sqrt{1+b^2+2b \sin(h s)}.
$$
It is clearly that $\Sig_b$ is invariant under the translations $z\mapsto z+x_b(2k\pi/h) e_1$, with  $e_1=(1,0,0)$ and $k\in \Z$.  
Moreover $\Sigma_0=\R\times \frac{1}{h} S^1 $, the straight cylinder of width $\frac{1}{h}$.   Furthermore, with the change of variables $ \cos(\frac{\th}{2})=\frac{h}{2} y_1(s)$, for $s\in(-\frac{\pi}{2h},\frac{\pi}{2h})$ and $\th=\th(s)\in (0,{\pi})$, we see  that 
$$
x_1(s)=\frac{\sqrt{2}}{h}-\frac{2}{h}\sin(\th) \qquad \textrm{ and }  \qquad y_1(s)=\frac{2}{h}\cos(\th).
$$
That is 
$$
 \Sig_1=\left\{ \frac{2}{h} \left(-\sin(\th), \cos(\th) \s\right)+\frac{\sqrt{2}+ 4k}{h} e_1\,:\, \th\in [0,\pi],\, \s\in S^1,\, k\in \Z \right\},
$$ 
 where $e_1=(1,0,0)\in \R^3$.  Therefore $\Sig_1$ is a family of tangent spheres with radius $\frac{2}{h}$, centered at the points $ \frac{\sqrt{2}+ 4k}{h} e_1$, $k\in \Z.$  \\
Next,  we observe  that for $b\in (0,1)$, the map $x_b$ is increasing on $\R$.  
  Hence with  the change of variable $t=x_b(s)$, and setting $\vp(t)=y_b(x_b^{-1}(t))$, we deduce that 
$$
 \Sig_b=\{ (t,\zeta)\in \R\times \R^2\,:\, t\in \R,\,|\zeta|= \vp_b(t) \}.
$$
Recently, Sicbaldi and Schlenk used in  \cite{ScSi} the Crandall-Rabinowitz   bifurcation theorem (see \cite{CR}) to derive these surfaces for $b$ close to $0$. We record it here. 
 \begin{Theorem}\label{Sicbaldi-Schlenk} There  exist    $b_0, h_*>0$ and a smooth curve $(-b_0,b_0)\ni b\mapsto \l(b)$    such that   $\l(0)=1$ and 
$$
 \vp_b (t)=\frac{1}{h_*}+\frac{b}{\l(b)}\left\{\cos\left( \lambda(b) t\right)+ v_{a}(\lambda(b)t)\right\},
$$
where $v_{b}\to 0$ in $C^{2,\a}(\R/{2\pi\Z})$ as $b\to 0$ and $\int_{-\pi}^{\pi}v_{b}(t)\cos(t)\, dt=0$ for every $b \in (-b_0,b_0)$.   
\end{Theorem}
 \begin{Remark}
We note that the family of surfaces   $(\Sig_b)_{b>1}$ are the  immersed constant mean curvature surfaces known as the \textit{nodoids}.
 \end{Remark}

Expression \eqref{eq:NMC-surf-rev} was first derived by  M. Sturm, in an appended note in\cite{Delaunay},  and    characterizes the Delauney
surfaces variationally,  as the extremals of surfaces of rotation having fixed volume
while maximizing lateral area.  A gereralization of Delauney surfaces in higher dimension is due to Hsiang and Yu, \cite{Hsiang-Yu}.

%
\section{Constant  Nonlocal Mean Curvature hypersurfaces}
Similarly to the mean curvature, the fractional or  Nonlocal Mean Curvature (NMC for short) is an extrinsic geometric quantity that is  invariant under global immersion  representing a surface.  Let $\alpha \in (0,1)$. If $\Sig$ is a smooth oriented hypersurface in $\R^N$ with unit normal vector field $\cV_\Sig$, 
its nonlocal   mean curvature  of order $\alpha$ at a point $x \in \Sig$
is defined as  
\begin{equation}
  \label{eq:def-frac-curvature}
H_\a(\Sig; x)=\frac{2}{\a} \int_{\Sig} \frac{(y-x)\cdot \cV_\Sig (y)}{|y-x|^{N+\alpha}}\,dV(y)  .
\end{equation}
 If $\Sig$ is of class $C^{1,\beta}$ for some $\beta>\alpha$ 
and we assume $\int_{\Sig} (1+|y|)^{1-N-\alpha}\,dV(y)< \infty$, then the integral in \eqref{eq:def-frac-curvature} 
is absolutely convergent in the Lebesgue sense. The orientation is chosen here so that $H_\a$ is positive for a sphere.  \\
We note that if $\Sig$ is of class $C^{2}$, then the normalized  nonlocal mean curvature 
 $\frac{1-\a}{  \o_N'  } H_\a(\Sig;\cdot)$  converges,  as $\a\to 1$, locally uniformly to the classical mean curvature $H(\Sig;\cdot)$ defined in \ref{eq:def-MC-geom}, see \cite{Davila2014B, Abatangelo}. Here $\o_N'$ is the measure of the $(N-2)$-dimensional unit sphere.\\
There is  an alternative expression for   $H_\a(\Sig;\cdot)$  in terms of  a solid integral.  Suppose  that $\Sig=\de E$ for some open set $E\subset\R^N$ and $\cV_{\de E}$ is the normal  exterior to $E$. Then, for all $x\in \Sig$,  we have 
\be  \label{eq:NMC-PV}
H_\a(\de E;x)=   PV\int_{\R^N}\frac{\t_{E}(y)}{|y-x|^{N+\a}}\, dy=  \lim_{\e\to 0 }\int_{\R^N\setminus B_\e(x)}\frac{\t_{E}(y)}{|y-x|^{N+\a}}\, dy,
\ee
where $\t_E:=1_{\R^N\setminus {E}}-1_{E}$, with    $1_D$ denotes the characteristic function of a set $D \subset \R^N$. This can be derived using the divergence theorem and the fact that  
\be \label{eq:div--y-x-weight}
\nabla_y\cdot (y-x){|y-x|^{-N-\alpha}} =-  {\alpha}{|y-x|^{-N-\alpha}}.
\ee

\begin{Remark}
The formula of NMC   in \eqref{eq:NMC-PV}  is comparable with  the one of the mean curvature, when written in infinitesimal solid integral,
$$
H(\de E;x): =   C  \lim_{\e \to 0}\frac{-1}{\e |B_\e(x)|} \int_{   B_\e(x)}\! \t_E(y) \,dy,
$$
for some constant  $C>0$ depending only on  $N$. Here the orientation is chosen so that the mean curvature is positive for a sphere. 
 In particular the mean curvature of $\de E$ is an infinitesimal  average, centered at $\de E$, of the sum of the "-1" coming from inside $E$   and the "+1" from outside $E$.  After all by Young's law, the mean curvature measures pressure difference across the   interface  of   two non mixing fluids at rest.
On the other hand the nonlocal mean curvature is a  weighted average of the sum of all the  "-1" coming from inside $E$  and all  the "+1" from outside $E$.
\end{Remark}

Both forms of the NMC \eqref{eq:def-frac-curvature} and  \eqref{eq:NMC-PV} turn out to be useful depending on the users interests.
 For instance,  global comparison principles is easily proved using \eqref{eq:NMC-PV}, while expression \eqref{eq:def-frac-curvature}  (without principle value integration) seems more convenient to work with when dealing with  regularity of  the NMC operator acting on graphs. By noticing that 
      if $E_1 \subset E_2$ then $\t_{E_1}- \t_{E_2}=2 1_{E_1\setminus E_2}$ on $\R^N$, we can state the following result. 
\begin{Lemma} \label{lem:glob-comp}
Let $E_1, E_2$ be two open sets of class $C^{1,\b}$, $\b>\a$, in a neighborhood of  $p\in \de E_1\cap \de E_2$. If $E_1\subset E_2$, then  $ H_\a(\de E_2;p)\leq H_\a(\de E_1;p)$, with equality if and only if $E_1=E_2$, up to set of  zero Lebesgue measure. 
\end{Lemma}      
 Local comparison principle   does not hold true in general as in the classical case.\\
 
Alike the mean curvature, the nonlocal mean curvature  appears  in the first variation of  nonlocal perimeter functional as well. The fractional perimeter of a bounded open set $E\subset\R^N$ is given by 
\begin{align*}
 P_\a(E)&=\frac{1}{2}  \int_{\R^N} \int_{\R^N }\frac{|1_E(x)-1_E(y)|^2 }{|y-x|^{N+\a}}\, dx dy=    \int_E \int_{\R^N\setminus E}\frac{1 }{|y-x|^{N+\a}}\, dx dy .
\end{align*}
We note that, if $E$ is a set with finite perimeter, the  normalized fractional perimeter $\frac{1-\a}{  \o_{N-1}  } P_{\a}(E)$  converges, as $\a\to 1$,  to $P(E)$, see \cite{Ambrosio,  Davila2014B}. Here $\o_{N-1}$ is the measure of  the unit ball in $\R^{N-1} $.\\
Isoperimetric type inequalities with respect to this functional was investigated in \cite{GR,FS}. 
According to a result of Frank and Seiringer \cite{FS}, it is known that balls uniquely minimize  $P_\a$ among all sets with a equal volume. A quantitative stability of the fractional isoperimetric inequality is proven by Fusco, Millot and Morini, \cite{FMM}.\\
Provided $E$ is bounded with Lipschitz boundary, using integration by parts and \eqref{eq:div--y-x-weight}, we can rewrite the fractional perimeter as 
\begin{align}
 P_\a(E)&=  \frac{ 1 }{\a(1+\a)} \int_{\de E} \int_{\de E}\frac{ \cV_{\de E}(y)\cdot( y-x)  \cV_{\de E}(x)\cdot( x-y)  }{|y-x|^{N+\a}}\, dV(y)   dV(x). \label{eq:second-express-of-frac-Per}
\end{align}
Note that \eqref{eq:second-express-of-frac-Per} provides a natural way to define a nonlocal or fractional measure of an orientable compact hypersurface.\\
%
In the theory of  minimal surfaces, the Plateau's problem is to show the existence of a surface that locally minimize perimeter with a given boundary.  For a nonlocal setting of  Plateau's problem,  in  \cite{Caffarelli2010} the authors introduced    the fractional perimeter in a reference open set $\O$  given by
\begin{align*}
 P_{\a,\O}(E)&=\frac{1}{2}  \iint_{\R^{2N}\setminus (\O^c)^2}\frac{ |1_E(x)-1_E(y)|^2}{|x-y|^{N+\a}}\, dxdy\\ &=L(E\cap\O,E^c\cap\O)+L(E\cap\O,E\cap\O^c)+L(E^c\cap\O,E^c\cap\O^c),
\end{align*}
where $A^c:=\R^N\setminus A$ and the interaction functional is given by 
$$
L(A,B):= \int_{B}\int_A\frac{dxdy}{|x-y|^{N+\a}}.
$$
We note   that the interaction between $E^c\cap \O^c  $ and $E\cap \O^c $ is left free. This, allows  a wellposed setting of  a (nonparametric) nonlocal Plateau's problem: existence  of  sets minimizing   $P_{\a,\O}$ among all sets that coincide in $\O^c$.  
The seminal paper \cite{Caffarelli2010} established the first
existence and regularity results for nonlocal perimeter minimizing sets in a  reference set $\O$, and   moreover  the boundary of the minimizers have,  in a viscosity sense, zero NMC in $\O$.   In the literature the boundary of such sets are called   \textit{$s$-minimal or nonlocal minimal hypersurfaces}.\\


The following formula for the  first variation of fractional perimeter was found  in \cite{FFMMM,Caffarelli2010 }.  
\begin{Theorem}  Let $\O$ and  $E$ be two domains of $\R^N$, with $E$ of class  $C^{1,\b}$, $\b>\a$. Let  $\l\in \R$ and $(E_t)_{t\in \R}$ be a variation of $E$ with respect to $\zeta\in C^\infty_c(\O;\R^N)$. Then the map $t\mapsto J_\a(t):= P_{\a,\O}(E_t)-\l |E_t\cap\O|$ is differentiable at zero. 
Moreover 
$$
J_\a'(0) =(N-1)\int_{\de E  } \left\{H_\a(\de E; p)- \l  \right\} v(p)\, dV(p) ,
$$
where   $v(p):= \la \zeta(p),  {\cV}_{\de E}(p)   \ra$.
 \end{Theorem}
The paper  \cite{FFMMM} contains also the second variation of fractional perimeter. \\
 From now on, we say that $\Sig$ is a CNMC hypersurface if it is of class $C^{1,\b}$, for some $\b>\a$, and such that $H_\a(\Sig,\cdot)$ is a \textit{nonzero} constant on $\Sig$.\\
Even though  we will be mainly interested in this note to CNMC hyperusfrcaes, we find it important to make a  brief digression in the  theory of the   nonlocal minimal hypersurfaces and those with vanishing NMC but not necessarily minimizing fractional perimeter, since it is also a hot topic nowadays with challenging open questions,  see Bucur  and Valdinoci \cite{Bucur}.\\
We first  recall that besides the hyperplane, which   trivially has zero  nonlocal mean curvature, there are some nontrivial  ones:      the $s$-Lawson cones found by D\'avila, del Pino and Wei \cite{Davila2014B}; the helicoid  found by Cinti, del Pino and D\'avila in \cite{CDD}. Moreover, local inversion arguments have been used in \cite{Davila2014B}, to derive,    for $\a$ close to 1,  two interesting examples of   rotationally symmetric surfaces with zero nonlocal mean curvature. The first one   posses the shape of the catenoid whereas the other is disconnected with two ends which are  asymptotic to a cone of revolution.  \\
Regularity theory, stability and  Bernstein-type results in the nonlocal setting are also studied in the recent years. Recall that  in the theory of constant mean curvature surfaces,   the boundary of perimeter     minimizing regions  are smooth except  a closed singular set of Hausdorff dimension at most $N-8$.   Such regularity result, in its full generality,  is still not  known to be true in the nonlocal setting. The  progress  made in this direction  so far parallel the classical regularity theory    up to a     \textit{possible dimension shift}. This latter fact    was discovered by D\'avila, del Pino and Wei  in \cite{Davila2014B} proving,  for $\a$  close to 0  and $N=7$,  that there are nonlocal stable minimal surfaces  ($s$-Lawson) cones --- but their globally minimizing property  is an open problem. On the other hand when $\a$ is close to 1,  Caffarelli and Valdinoci~\cite{Caffarelli2011A,Caffarelli2011B}
proved the $C^{1,\gamma}$ regularity of nonlocal minimal hypersurfaces, except a set of $(N-8)$-Hausdorff dimension.  This, together with the subsequent results of Barrios, Figalli, and Valdinoci~\cite{Barrios}, 
leads, for $\alpha$ close to 1, to their $C^{\infty}$ regularity up to dimension $N\leq 7$. We note that for $N=2$ and any $\alpha\in (0,1)$, Savin and Valdinoci~\cite{Savin2012} established  that nonlocal minimal curves  are
$C^{\infty}$. Decisive regularity estimates have been  proved  in \cite{Caffarelli2010, Figalli2013, CSV,CCS,CC }, see also 
the monograph \cite{Bucur} and the references there in.\\

We close this section by noting that nonlocal capillary problems have been investigated in \cite{MV,Mihaila} and nonlocal mean curvature flow in \cite{Saez2015, Cinti.Sinestrari.Valdinoci-ArXiv,Imbert2009}. 

 \subsection{Expressions of  the  NMC of some globally parameterized hypersurfaces  }
 In this section, we derive some formulae of some hypersurfaces which admit    global parameterizations.  We will use both expression of the NMC in \eqref{eq:def-MC-geom} and \eqref{eq:def-frac-curvature}. Depending on the problems under study, each expression has its own advantages. For instance to prove regularity of the NMC operator, it is in general more convenient to consider those formulae  from  \eqref{eq:def-MC-geom}, whereas those  from \eqref{eq:def-frac-curvature} turns out to be useful in the study of qualitative properties of the  NMC acting on graphs as  a quasilinear nonlocal elliptic operator.    
 \subsubsection{Without principle value integration}
We have the following result.
\begin{Proposition} \label{prop:NMC-gen-slabs}  
Let $N=n+m\geq 1$, with  $n,m\in \N$. 
Let ${u}:\R^m\to (0,+\infty)$ be  a function of class $C^{1,\b}$, $\b>\a$,  and  satisfy 
$$
\int_{\R^m} \frac{(1+|\n u (\t)|)u^{n-1}(\t)}{(1+ |\t|+ u(\t))^{N-1+\a}}d\t<\infty.
$$ 
  Consider    the set 
$$
E_u:=\{ (s,\z)\in \R^m\times \R^{n}\,:\, |\z|<u(s)\}.
$$
\begin{itemize}
\item[(i)] For $n\geq 2$,   the NMC of $\Sig_u$ at the point $q=(s,u(s) e_1)$ is   given by 
\begin{align}
  -\frac{\alpha}{2}  H_\a(\Sig_u;q ) =& \int_{S^{n-1}}  \int_{\R^m}
  \frac{
\bigl\{ u(s)-u(s-\tau)-\tau \cdot  \n u(s-\tau) \bigr\}u^{n-1} (s-\t)
}{
\{|\t|^2+({u}(s)-{u}(s-\t))^2+{u}(s){u}(s-\t)|\sigma-e_1|^2 \}^{(N+\a)/2}
}
d\tau d\s   \nonumber \\
& \hspace{-14mm} -
 \frac{{u}(s)}{2} \int_{S^{n-1}}   \int_{\R^m}
 \frac{
|\sigma-e_1|^2 u^{n-1} (s-\t)
}{
\{|\t|^2+({u}(s)-{u}(s-\t))^2+{u}(s){u}(s-\t)|\sigma-e_1|^2 \}^{(N+\a)/2}
}
d\t d\s  , \label{new-rep-H}
\end{align} 
where $e_1=(1,0,\dots,0) \in \R^{n}$.
\item[(ii)] For $n=1$, the NMC of $\Sig_u$ at the point $q=(s,u(s))$ is given by 
\begin{align}
  -\frac{\alpha}{2}  H(\Sig_u;q) =&    \int_{\R^{N-1} }
  \frac{
  u(s)-u(s-\t)-\t \cdot  \n u(s-\t)   
}{
\{|\t|^2+({u}(s)-{u}(s-\t))^2  \}^{(N+\a)/2}
}
d \t\label{eq:Geom-from-NMC} \\
& -   \int_{\R^{N-1} }
  \frac{
  u(s)+u(s-\t)+  \t \cdot  \n u(s-\t) 
}{
\{|\t|^2+({u}(s)+{u}(s-\t))^2  \}^{(N+\a)/2}
}
d \t  .  \nonumber   
\end{align} 
\end{itemize}
Moreover,   all integrals above converge absolutely in the Lebesgue sense.
\end{Proposition}
Expression \eqref{new-rep-H}  and \eqref{eq:Geom-from-NMC} are easily  derived from \eqref{eq:def-frac-curvature},  by  change of variables, taking into account that  the unit outer normal of   $\de \Sig_{u}$ at the point $q=(s,u(s)\s ) $ is given by
 $$
\cV_{\de E_{u}}( q)=\frac{1}{\sqrt{1+|{\n u}|^2({s})}}(-{\n  u}(s),\s) 
 $$
and the volume element is ${u}^{n-1}({s}) \sqrt{1+|{\n u}|^2({s})}dsd\s .$ The proof uses similar arguments as in \cite{CFW-2016}. We note that the first term in the right hand of   \eqref{eq:Geom-from-NMC} provides  the expression of NMC of an euclidean   graph $x_N=u(x')$ without principle value integral.
%

  \subsubsection{With  principle value integration}
  We consider a function $u:\R^{N-1}\to \R$ of class $C^{1,\b}$, for some $\b>\a$.   We let
$$E_u=\left\{(y,t)\in\R^{N-1}\times\R\,:\, t  < u(y) \right\}$$
and consider the  parametrization   
$$
\R^N\to \R^N\qquad (y,t)\mapsto \calF(y,t)=\left( y,u(y)+t \right).
$$
It is well known that this parametrization is \textit{volume preserving}.  
We now compute the  NMC of the graph $u$ (denoted by $H(u)$) at the point $\calF(x,0)=(x,u(x))\in\de E_u$.  From  \eqref{eq:def-frac-curvature}, making change of variables and using Fubini's theorem, we have 
\begin{align*}
H(u)(x):=H_\a(\de E_u;(x,u(x))) 
 &=\int_{\R^{N}}\frac{\tau_{E_0}(y,t)}{|(x,u(x))-\calF(y,t)|^{N+\a}}  d ydt\\
  &=\int_{\R^{N-1}}I(x,{y}) dy,
\end{align*}
where 
\begin{align*}
I({x},{y})&=[\int_{-\infty}^{0}- \int^{+\infty}_{0}]\left\{  |{x}-{y}|^2+(t+u({y})-u({x}))^2 \right\}^{\frac{-(N+\a)}{2}}dt\\
&=|{x}-{y}|^{-(N+\a)} [\int_{-\infty}^{0}- \int^{+\infty}_{0}]
\left\{ 1  + \left(\frac{t+u({y})-u({x})}{ |{x}-{y}|}\right)^2 \right\}^{\frac{-(N+\a)}{2}}d t.
\end{align*}
 We then make the change of variables $s=\frac{t+u({y})-u({x})}{ |{x}-{y}|}$ and
define 
 \be \label{eq:P_u}
  p_u({x},{y})= \frac{u({y})-u({x})}{|{x}-{y}|} ,
 \ee
to get
\begin{align*}
I({x},{y}) &=-|{x}-{y}|^{-(N-1+\a)} [\int_{-\infty}^{p_u({x},{y})}- \int^{+\infty}_{p_u({x},{y})}]
\left\{ 1  + s^2 \right\}^{\frac{-(N+\a)}{2}}ds.
\end{align*}
Lgtting
\be \label{eq:def-of-F}
F(p):=\int_p^{+\infty}\frac{d\t}{(1+\t^2)^{\frac{(N+\a)}{2}}},
\ee
we find that 
$$
 H(u)({x})= \int_{\R^{N-1}}\frac{F(p_u({x},{y}))- F(-p_u({x},{y})) }{|{x}-{y}|^{N-1+\a}} d{y}.
$$
We then have the following formulae for  the NMC of the hypersurface  $\de E_u$.
\begin{Proposition}\label{prop:exp-NMC-pv1}
Let $u: \R^{N-1}\to \R$ of class $C^{1,\b}$, $\b>\a$. Consider the subgraph of $u$,  
$$
E_u=\left\{(y,t)\in\R^{N-1}\times\R\,:\, t  < u(y) \right\} .
$$
At a point  $q=(x,u(x))\in \de E_u$, we have 
\begin{align}
 H(u)(x):=H_\a(\de E_u;q)&=PV\int_{\R^{N-1}}\frac{F(p_u({x},{y}))- F(-p_u({x},{y})) }{|{x}-{y}|^{N-1+\a}} d{y} \label{eq:NMC_curve-E1}\\
&=   PV\int_{\R^{N-1}}\frac{u({x})-u({y})}{|{x}-{y}|^{N+\a}}\,q_u(x,y) \, d {y},
 \label{eq:NMC_curve-E2}
   \end{align} 
   where  $p_u$ and $F$ are  given by \eqref{eq:P_u} and  \eqref{eq:def-of-F}, respectively, and  
$$
q_u(x,y)
=\int_{-1}^1 {\left( 1+\t^2\frac{(u(x)-u(y))^2}{|x-y|^2} \right)^{-\frac{N+\a}{2}}}d\t.
$$
  \end{Proposition} 
  Here, \eqref{eq:NMC_curve-E2},  is a consequence of   the fundamental theorem of calculus, which yields
$$
F(p)-F(-p)= p\int_{-1}^1F'(\t p)\, d\t.
$$
As in the classical case, we may expect that the difference of  the NMC operator of  two graphs gives rise to a  nonlocal symmetric operator   allowing for  local comparison principles.  
\begin{Lemma}\label{lem:local-comparison-NMC-graphs}
Let $\O$  be an open set of  $\R^{N-1}$ and   $u,v\in C^{1,\b}_{loc}(\O)\cap    C(\R^{N-1})$. Then letting $w=u-v$, for every $x\in \O$, we have
 \be\label{eq:Hu--Hv}
   H(u)(x)-H(v)(x)=PV\int_{\R^{N-1}}\frac{w(x)-w(y)}{{|x-y|^{N+\a}}} \ti q_{u,v}(x,y)  \,dy,
 \ee
where $\ti q_{u,v}(x,y):=-2  \int_{0}^{1}F'\left(p_v(x,y)+\rho p_{u-v}(x,y)\right)\, d\rho.$
In particular,  if 
$$
 H(u)(x)-H(v)(x)\geq 0\qquad \textrm{ for every } x\in \O
$$
and $u\geq v$ on $\R^{N-1}\setminus \O $,  then $u-v$ cannot  attain its global minimum in $\O$ unless it is constant on $\R^{N-1}$.
\end{Lemma}
\proof
 Since    $F'$ is even,     by  the fundamental theorem of calculus, we get
$$
[F(a)-F(b)]- [F(-a)-F(-b)]=2(a-b)\int_0^1 F'\left(b+\rho(a-b)\right) d\rho.
$$
In view of \eqref{eq:NMC_curve-E1},   this gives \eqref{eq:Hu--Hv} since  $ p_w({x},{y})=- \frac{w({x})-w({y})}{|{x}-{y}|} $.\\
For the comparison principle, we  suppose that, for some $x_0\in\O$, we have  $w(x_0)= \min_{x\in \R^ {N-1}}w(x).$
We then have  
\begin{align*}
 0&\leq   H(u)(x_0)-H(v)(x_0)\\
 & =-2 PV \int_{\R^{N-1}}\frac{w(x_0)-w(y)}{{|x_0-y|^{N-1+\a}}} \int_{0}^{1}F'\left(p_v(x_0,y)+\rho p_w(x_0,y)\right)\, d\rho  \,dy \leq  0.
\end{align*}
Therefore, since $F'<0$  on $\R$,   we deduce that $w\equiv w(x_0)$ on $\R^{N-1}$.
\QED


%
%
In the case of   generalized slabs, we have the following result, the proof is similar to one in \cite{CFW-2016}.
\begin{Proposition}\label{prop:exp-NMC-pv2}
Let ${u}:\R^m\to (0,+\infty)$ be a a function  of class $C^{1,\b}$ and $N=m+n$.  Consider the open set
$$
E_u:=\{ (s,\z)\in \R^m\times \R^{n}\,:\, |\z|<u(s)\}.
$$
Then at the point  $q=({s}, u({s}) z)$, we have 
\begin{align*}
H_\a(\de E_u;q) 
&= PV \int_{\R^n} \int_{\R^{m}}\frac{-\t_{E_1}(\bar{s},z)}{((s-\bar{s})^2+ (\u(s)e_1-u(\bar{s})z)^2  )^{\frac{N+\a}{2}}}u^{n}(\bar{s})\,dz d\bar{s} .
\end{align*}
\end{Proposition}
  We note that in the expressions of the NMC with PV in Proposition \ref{prop:exp-NMC-pv1} and \ref{prop:exp-NMC-pv2}, no growth control of $u$ at infinity is required.

 \subsection{Bounded constant nonlocal mean curvature hypersurfaces}
In addition to the cylinders, spheres, which are CNMC with nozero NMC,  we shall see that there many more. However,    this class  is   reduced by    a nonlocal  counterpart of the   Alexandrov  result 
on the characterization of spheres as the only closed embedded CMC-hypersurfaces. 
 
\begin{Theorem}[\cite{Cabre-Alex-Del-2015A, Ciraolo2015}]
\label{th:nonlocal-Alex}
Suppose that $E$ is a nonempty bounded open set (not necessarily connected) with $C^{2,\beta}$-boundary for some $\beta >\alpha$ and with the property that $H_\a(\de E;\cdot)$ is constant on $\partial E$. Then $E$ is a ball. 
\end{Theorem}
This result was obtained at the same time and independently by Cabr\'e, Sol\`a-Morales, Weth and  the author in \cite{Cabre-Alex-Del-2015A} and by Ciraolo, Figalli, Maggi, and Novaga \cite{Ciraolo2015}.
We mention that the  paper  \cite{Ciraolo2015} contains also 
stability results with respect to this rigidity theorem.\\
This new characterization of the sphere  relies on the Alexandrov's moving planes method   discussed above.   Indeed,    pick a unit vector $e\in \R^N$ together with a hyperplane $P_e$ which is perpendicular to $e$ and not intersecting $E$. Call $E_e$ the reflection of $E$ with respect to the plane $P_e$. 
Then slide the plane along   the $e$-direction, toward $E$, until one of the two properties occurs for the first time:
\begin{enumerate}
\item[1.] \textit{interior touching:} there exists $x_0\not \in P_e$, with   $x_0\in \de E\cap \de  E_e$,
\item[2.] \textit{non-transversal intersection:}  $e\in T_{x_0} \de E $, for some $x_0\in \de E$.
\end{enumerate}
In both cases, we will  find that $E_e=E$, and since $e$ is arbitrary, this implies that $E$ must be a unit ball.
As mentioned earlier,   in contrast to the classical case,  there is no local comparison principle related to the fractional mean curvature. Moreover, while comparison principle holds for graphs (see Lemma \ref{lem:local-comparison-NMC-graphs}), we dare not hope for  a global parameterization of $E$ by a graph!   However,  we may rely on  a   global comparison principles inherent to problem. The montonicity of the weight $t\mapsto K_N(t):= |t|^{-N-\a}$ will be crucial.  Indeed,  in  case \textit{1.},  we have 
\begin{align*}
0&= H_\a(\de E;x_0)-H_\a(\de E_e;x_0)=\frac{1}{2}PV\int_{\R^N}( \t_E(y)-\t_{E_e}(y))K_N(x_0-y)  \,dy\\
&=\frac{1}{2}PV\int_{E\setminus E_e}( \t_E(y)-\t_{E_e}(y))K_N(x_0-y)\,dy+ \frac{1}{2}PV\int_{E_e\setminus E}( \t_E(y)-\t_{E_e}(y))K_N(x_0-y)\,dy\\
&=\frac{1}{2}PV\int_{E\setminus E_e}K_N(x_0-y)\,dy- \frac{1}{2}PV\int_{E_e\setminus E} K_N(x_0-y)\,dy\\
&=\frac{1}{2}PV\int_{E\setminus E_e}\left( K_N(x_0-z)-   K_N(x_0-\cR(z))\right)\,dz.
\end{align*}
We made the  change  of variable,     $y= \cR_e(z)$,  with $\cR_e$ being the reflection with respect to the plane  $P_e$. 
Now since $ E\setminus E_e$ is  contained on the side of the plane containing $x_0$, we find that $ |x_0-z| <   |x_0-\cR(z)|  $ for every $z \in E\setminus E_e$,  together with the monotonicity of $K_N$  imply that  $E=E_e$.  \\
Case \textit{2.} is far more tricky  and requires  a formula for directional   derivative  of $H_\a(\de E;\cdot)$, see
Proposition~\ref{sec:introduction-1} below. Since $e\in T_{x_0}\de E=T_{x_0}\de E_e$,  we have 
\begin{align*}
0 &= \frac{1}{2}\bigl(\partial_{e} H_\a(\de E;x_0)-\partial_{e} H_\a(\de E_e;x_0)  \bigr)\\
 &= -\frac{N+\a}{2}  PV \int_{\R^N} (x_0\cdot e-y\cdot e) \bigl(\tau_E(y)- \tau_{ E_e}(y)\bigr)  K_{N+2}(x_0-y) \,dy\\
 &= -\frac{N+\a}{2}PV\int_{\R^N}(x_0\cdot e-y\cdot e) (1_{E\setminus E_e }(y) - 1_{E_e\setminus E }(y)  )K_{N+2}(x_0-y)\,dy .
\end{align*}
Since the integrand does not change sign on $\R^N$,    there must be $E=E_e$. We then conclude that in either case   $E=E_e$, and since $e$ is arbitrary, we deduce that  $E$ is a ball.\\
Of course to make all these argument rigorous, one should take into account that the integrals are defined in the principle value sense, see \cite{Cabre-Alex-Del-2015A,Ciraolo2015} for more details.

The following   formula for the derivatives of the nonlocal mean curvature was used above, see \cite{Cabre-Alex-Del-2015A,Ciraolo2015}.
\begin{Proposition} 
\label{sec:introduction-1}
If $E \subset \R^N$ is bounded and $\partial E$ is of class $C^{2,\beta}$ for some $\beta> \alpha$, then $H_E$ is of class $C^1$ on $\partial E$, and we have 
$$
\partial_v H_\a(\de E;x) =- (N+\alpha)\, PV \int_{\R^N} \tauE(y)  |x-y|^{-(N+2+\alpha)} (x-y) \cdot v \,dy
$$
for $x \in \partial E$ and $v \in T_x \partial E$.
\end{Proposition}
 
In order  to reduce to a single sphere,   connectedness   is obviously a necessary assumption in Theorem~\ref{th:Alexandrov} whereas Theorem \ref{th:nonlocal-Alex} allow for disconnected  sets $E$ with finitely many connected components.  Similar pheonomenon was also observed in the study of fractional  overdetermined problems, see  \cite{FJ}.

 \subsection{Unbounded constant nonlocal mean curvature hypersurfaces.} 
 When compactness is dropped, the cylinder provides a trivial  example of surfaces with constant and nonzero local/nonlcal mean curvature.      Obviously,   CMC curves bounding a periodic domain    in $2$-dimension are  parallel straight lines and hence has zero boundary mean curvature. This latter fact does not hold in the nonlocal case, since two parallel straight lines in $\R^2$ has      positive CNMC. In fact, any slab $\{(s,\z)\in \R^n\times \R^m\,:\, |\z|=1\}$ has a positive CNMC, as it can be easily seen  from  Proposition \ref{prop:NMC-gen-slabs}.  It is therefore  natural to  ask   if there are CNMC hypersurfaces besides the cylinder?\\
The lack of ODE theory in the nonlocal framework make this question not trivial to answer. Up to now, all existence  results use either  variational methods or  perturbation methods.  
 
   The first answer to the above   question, for $N=2$,  has been given by  Cabr\'{e},   J. Sol\`{a}-Morales,   Weth and the author  in    \cite{Cabre-Alex-Del-2015A}, where   a continuous branch of periodic CNMC sets bifurcating from the straight band was found.  This was improved and  generalized in \cite{CFW-2016} to higher dimensions.   
We note that  in \cite{Davila2015},   D\'avila, del Pino, Dipierro, and Valdinoci, established  variationally
the existence of periodic and cylindrically symmetric hypersurfaces in $\R^N$
which minimize the periodic  fractional perimeter under a volume constraint.
More precisely, \cite{Davila2015}~establishes the existence of a $1$-periodic minimizer for every given volume within the
slab $\{(s,\z)\in\R\times\R^{N-1}\, : \, -1/2<s<1/2\}$, which are CNMC  hypersurfaces in  the viscosity  sense of \cite{Caffarelli2010}, since  it is not known if they are of class    $C^{1,\beta}$ for some $\beta>\alpha$.
\subsubsection{ CNMC hypersurface  of revolution }
We consider hypersurfaces   with constant nonlocal mean curvature
of the  form 
\begin{equation*}
  \label{eq:cylindrical-graph}
\Sig_u=\{(s,\z)\in \R\times \R^{N-1}\, : \, |\z|=u (s)\},
\end{equation*}
where $u: \R \to (0,\infty)$ is a positive and even function. 
The following result shows  the existence of a smooth branch of sets  
which are periodic in the variable~$s$ and have
all the same constant nonlocal mean curvature; they bifurcate from a straight cylinder $\Sig_R:=\{|\z|=R\}$. 
The radius $R$ of the straight cylinder is chosen so that the periods of the new cylinders converge to $2\pi$ as they
approach the straight cylinder.   We state the counterpart of Theorem \ref{Sicbaldi-Schlenk}.
\begin{Theorem}
\label{res:cyl1}
Let $N \ge 2$. For every $\alpha\in(0,1)$, $\b\in (\a,1)$ there exist $R, a_0>0$ and  smooth curves $a\mapsto u_a$ and $a\mapsto \l(a)$, with $\l(0)=1$,  such that,   for every $a\in (-a_0,a_0)$,   
$$
 \Sig_{u_a}=\{(s,\z)\in \R\times \R^{N-1}\, : \, |\z|=u_a (s)\}
$$ 
 is a CNMC hypersurface of class $C^{1,\b}$, with 
$$
 H_\a( \Sig_{u_a};q )=  {H_\a( \Sig_{R};\cdot )} ,\qquad\textrm{ for all $q\in \de  \Sig_{u_a}$.}
$$
Moreover   for every $a \in (-a_0,a_0)$, we have 
\begin{equation*}\label{form}
  u_{a}(s)=R+\frac{a}{\l(a)}\left\{\cos\left( \lambda(a) s\right)+ v_{a}(\lambda(a)s)\right\}, \qquad u_{-a}(s)= u_a\Bigl(s+ \frac{\pi}{\lambda(a)}\Bigr),
\end{equation*}
where $v_{a}\to 0$ in $C^{1,\b}(\R/{2\pi\Z})$ as $a\to 0$ and $\int_{-\pi}^{\pi}v_{a}(t)\cos(t)\, dt=0$ for every $a \in (-a_0,a_0)$.
\end{Theorem}
The proof of Theorem \ref{res:cyl1} rests on the application the Crandall-Rabinowitz bifurcation theorem \cite{CR},  applied to    
the  quasilinear type fractional elliptic equation, see  Proposition \ref{prop:NMC-gen-slabs}, 
$$
 H_\a( \Sig_{u}; \cdot)-  {H_\a( \Sig_{R};\cdot )}=0 \qquad\textrm{  in  $  \de  \Sig_{u}$}.
$$
To do so and to obtain a smooth  branch of bifurcation parameter, we need  the smoothness of  $u\mapsto (s\mapsto  H_\a( \Sig_{u};(s,u(s) e_1)  ))$ as map from an open subset of  $C^{1,\b}(\R)$ taking values in $C^{0,\b-\a}(\R)$.  This is a nontrivial task, and as such the use of the formula for the NMC without PV (see Proposition \ref{prop:NMC-gen-slabs}) become crucial. \\

The first question that may come  to the reader's mind, is the existence of global continuous $\left( \Sig_{u_a}\right)_{ a\in (0,1)}$, 
branches of nonlocal Delaunay hypersurfaces and the transition from unduloids to a periodic array of spheres  as   in  the local case ($\a=1$) discussed in Section \ref{s:Del-local}.  It is indeed an open question.    However some key differences are expected in the embedded regime $a\in (0,1)$.  Indeed,  two results suggest that as the bifurcation parameter  $a$ varies from $0$ to $1$,  the hypersurfaces $\Sig_{u_a}$  should approach an infinite compound  of not-round-spheres.  The first one being  \cite{Davila2015}, where the authors proved     that the enclosed sets of  their  (weak)   CNMC Delauney hypersurface,     as their constraint volume goes  to zero,   
tend in measure (more precisely, in the so called Fraenkel asymmetry) to a
periodic array of balls. The second one is a consequence of the work  of Cabr\'e, Weth  and the author  in \cite{CFW-2017}, where it is proven that an array of periodic  disjoint round spheres is not necessary a CNMC hypersurface.  We detail  a bit more on the  latter result   in the   following section. 

\subsubsection{Near-sphere lattices with CNMC}
In this section, we consider  CNMC hypersurfaces given by an infinite compound of aligned round spheres, 
tangent or disconnected, enlightening possible limiting configuration of nonlocal unduloids $\Sig_{u_a}$ (from Theorem \ref{res:cyl1}) as the bifurcation parameter becomes large.    In a more general setting, we can look for  CNMC hypersurfaces which are   countable union of a certain bounded domain.
To be precise, we assume $N\geq 2$ and let 
$$
S:= S^{N-1}\subset \R^N
$$  
denote the unit centered sphere of $\R^N$.  For $M\in\N$ with $1\leq M\leq N$ we regard $\R^M$ as a subspace of $\R^N$ by 
identifying $x'\in \R^M$ with $(x',0)\in \R^M \times \R^{N-M} =  \R^N$.  Let 
$\left\{ \textbf{a}_1;\dots;\textbf{a}_M\right\}$  be a basis of~$\R^M$.   
 By the above identification, we then consider the $M$-dimensional lattice 
\be \label{eq:def-scrL-lattice}
\scrL=\left\{\sum_{i=1}^M k_i \textbf{a}_i\,:\,  k=(k_1,\dots,k_M)\in \Z^M\right\}
\ee 
as a subset of $\R^N$. In the case where     $\left\{ \textbf{a}_1;\dots;\textbf{a}_M\right\}$   is an   orthogonal or an orthonormal basis, we say that $\scrL$ is a  {rectangular lattice} or a square lattice, respectively. 
 We define, for $r>0, $  
$$
\scrS^r_0:=  S+ r\scrL:= \bigcup_{p \in\scrL} \Bigl( S+  r p \Bigr)  \; \subset \; \R^N.
$$
For $r$ large enough,  the set $\scrS^r_0$ is the  union of disjoint unit spheres  
centered at the  lattice points in $r \scrL$. Consequently,  $\scrS^r_0$  is a set of constant classical mean curvature 
(equal to one).  In contrast,   we shall see that
the NMC function $H_\a(\scrS^r_0;\cdot )$ is in general \textit{not} constant on $\scrS^r_0$. However nearby $\scrS^r_0$, one may expect CNMC sets for $r$ large enough. Indeed,  as $r\to \infty$, the hypersurface  $\scrS^r_0$ tends to the single centered sphere $S$, which has CNMC. On the other hand,  by invariance under   translations $\R^N$, the linearized NMC operator about $S$ has an $N$-dimensional kernel spanned   by the coordinates functions $x_1,\dots,x_N$. Taking advantages on the invariance of the NMC operator by even reflections, we shall see that   this program is realizable.    Indeed, we consider the open set 
$$
\cO:=\{\vp \in  C^{1,\b}(S) \,:\, \|\vp\|_{L^\infty(S)}<1, \textrm{   $\vp$ is even on $S$} \},
$$
with $\b\in (\a,1)$, and  the deformed sphere 
$
S_\phi:= \{ (1+ \phi(\s)) \s \::\: \s \in S\}, \qquad \vp \in \cO.
$
Provided that $r>0$ is large enough, the deformed sphere lattice  
$$
\scrS^r_\vp:=  S_\vp+ r\scrL := \bigcup_{p \in \scrL} \Bigl(S_\vp+  rp \Bigr)
$$
is a noncompact periodic hypersurface of class $C^{1,\beta}$. We have the following result. 
\begin{Theorem}[\cite{CFW-2017}]\label{th:main-th-Brav-lat}
Let $\a\in (0,1)$, $\b\in \left(\a,1\right)$,  $N\geq 2$, $1\leq M\leq N$ and $\scrL$ be an $M$-dimensional lattice. Then,
there exist $r_0>0$,  and a  $C^2$-curve $(r_0,+\infty) \to \cO$, $r \mapsto \phi_r$ such that 
 for every $r\in (r_0,+\infty)$,  the hypersurface  $S_{\vp_r}+ r\scrL $ has constant nonlocal mean curvature given by $H_\a(S_{\vp_r}+ r\scrL;\cdot  ) \equiv   H_\a(S;\cdot)$.\\
 Moreover if $1\leq M\leq N-1$,  then the functions $\vp_r$, $r > r_0$, are non-constant on $S$.
\end{Theorem}
As a consequence a periodic array of aligned spheres does not necessarily  have    constant nonlocal mean curvature.   A Taylor expansion of the perturbation  $\vp_r$, as $r\to \infty$,     shows how the  near-spheres interact to  form a CNMC hypersurface.  To see this, we consider  the linearized operator for the NMC operator acting on graphs  on  the unit sphere $S$, given by 
\begin{equation*}\label{eq:linNMC}
\phi \mapsto 2 (L_\alpha \phi-\lambda_1\phi),
\end{equation*}
where
\begin{equation*}
  \label{eq:def-L-alpha}
L_\a\vp(\th )= 
PV \int_S\frac{\vp(\th)-\vp(\s)}{|\th -\s|^{N+\a}}\, dV(\s).
\end{equation*}
The operator $L_\a$ can be seen as a \textit{spherical fractional  Laplacian}, and the above integral is understood in the principle value sense.
It  has the spherical harmonics   as eigenfunctions corresponding to an increasing sequence of eigenvalues $\l_0=0<\l_1<\l_2<...$.\\
The function $\vp_r$ in Theorem~\ref{th:main-th-Brav-lat} expands as 
\begin{align*}
\phi_r(\th) 
&=r^{-N-\alpha}\left(-\kappa_0  +   r^{-2}\left\{\kappa_1 \sum_{ p\in \scrL_*} \frac{(\th \cdot p)^2}{|p|^{N+\alpha+4}} -\kappa_2\right\}
 +  o\left(r^{-2}\right) \right) \quad \text{for $\theta \in S$ as $r \to +\infty,$}
\end{align*}
where  $\scrL_*:= \scrL \setminus \{0\}$, for some positive constants $\k_0,\k_2,\k_3$, only depending on $N,\a,\l_1,\l_2$ and on $\scrL$. Since $\kappa_0>0$, the above  expansion  shows
that, for large $r$, the perturbed spheres $S_{\vp_{r}}$ become smaller than $S$ 
as the perturbation parameter $r$ decreases.  
With regard to the order $r^{-N-\alpha}$, the shrinking process is uniform on $S$, whereas non-uniform deformations of the spheres may appear at the order $r^{-N-\alpha-2}$.   In the case $M=N$, it is not known whether $\vp_r$ is not constant on $S$. In fact,   we have the following more explicit  form of $\vp_r$ in   the case of  square lattices. 
Assume that  $\scrL$ is    a square  lattice.  Then
$$
\phi_r(\th )= r^{-N-\alpha}\left(-\kappa_0  +   r^{-2}\left\{ \tilde \kappa_1  \sum_{j=1}^M    \th_j^2 -\kappa_2 \right\}
 +  o(r^{-2}) \right) \quad \text{for $\th \in S$ as $r \to +\infty$},
$$
where $\tilde \kappa_1=\frac{\k_1}{M} \sum \limits_{p \in \scrL_*}\frac{1}{|p |^{N+\alpha+2}}$. In particular, if $M=N$, then 
$$
 \phi_r(\th )= r^{-N-\alpha}\Bigl(-\kappa_0  +   r^{-2}\Bigl(\tilde  \kappa_1  -\kappa_2 \Bigr)
 +  o(r^{-2}) \Bigr) \qquad \text{as $r\to \infty$.}
$$
Hence the deformation of the lattice $S_{\vp_r}+r\scrL$ is uniform up to the order $r^{-N-\alpha-2}$.   It is conjectured in \cite{CFW-2017} that   that $H_\a(\scrS^r_0;\cdot)$ is non-constant for any $N$-dimensional lattice $\scrL$, as long as $\scrS^r_0$ is a hypersurface of class $C^{1,\b}$, $\b\in (\a,1)$.\\
The proof of Theorem \ref{th:main-th-Brav-lat}, is based on the application of the implicit function  theorem, to solve the  $r$-parameter dependence    fractional quasilinear elliptic type problem on the sphere,
$$
H_\a( \scrS^r_\vp; \cdot )=  H_\a(S _\vp;\cdot ) + r^{-N-\a}G(r,\vp; \cdot)= H_\a(S _0;\cdot ).
$$
The   term containing $G$ in the above identity is  a  lower order term. In fact  the function $G$ together with all its derivatives in the $\vp$ variables are bounded, for $r$ large enough. On the other hand the leading quasilinear  operator   given by $ H_\a(S _\vp;\cdot )$ is 
 the NMC operator acting on graphs on the sphere. Its expression,  derived from \eqref{eq:def-frac-curvature},  is explicitly  given in the following %
\begin{Proposition}\label{prop:new-express-nmc}
Let $\psi\in C^{1,\b}(S^{N-1})$ with $\psi>0$. Then the NMC of the hypersurface
  $$\Sig_\psi:=\left\{\s\psi(\s)\,:\, \s\in S^{N-1} \right\}$$
  at the point $q=\th\psi(\th)$ is given by  
  \begin{align*} 
h_\a&(\psi)(\th):=\frac{\a}{2 } H_\a(\Sig_\psi  ; q)=-\psi(\th) \int_{{S}}\frac{  \psi(\th)-\psi(\s)   -  (\th-\s) \cdot\n\psi(\s) }{ |\th-\s|^{N+\alpha}}  \psi^{N-2}(\s)\,  \cK(\psi,\s,\th)       \,  dV(\s) \\
&+ \int_{{S}}\frac{(\psi(\th)-\psi(\s)) ^2  }{ |\th-\s|^{N+\alpha}} \psi^{N-2}(\s)\,  \cK (\psi,\s,\th)       \,  dV(\s) +\frac{\psi(\th)}{2}  \int_{{S}}\frac{ \psi^{N-1}(\s)  }{  |\th-\s|^{N+\alpha-2}}         \,   \cK (\psi,\s,\th)   \,dV(\s),\nonumber\\
\end{align*}
where
$
\cK(\psi,\s,\th):=\frac{1}{\left( \frac{(\psi(\th)-\psi(\s))^2}{|\th-\s|^2} +\psi(\s) \psi(\th)\right)^{(N+\a)/2}}.
$
Moreover, all integrals above converge absolutely.

\end{Proposition}
 Establishing the  regularity  of the NMC  operator $h_\a$, appearing in Proposition \ref{prop:new-express-nmc},  as a map from open subsets of  $C^{1,\b}(S^{N-1})$ and taking values in $ C^{0,\b-\a}(S^{N-1})$,   turns out to be an involved task. The   inconveniences in the above  expression  is the presence of  $\th$ in  the singular kernel $|\th-\s|^{-N-\a}$ and the fact that it involves only    euclidean distance instead of geodesic distance on the sphere.

\section{ Serrin's overdetermined problems}
We consider the problem of finding domains (not necessarily bounded) and functions  $u\in C^2(\ov \O)$ such that  
\begin{align}\label{eq:problemb}
    -\Delta_{g} u =1 \quad \textrm{ in }  \Omega, \qquad\qquad
u =0,   \quad \de_\nu u =const.\qquad\textrm{ on }\partial\Omega,
%
  \end{align}
   A domain $\O$  is called a \textit{Serrin domain} if it is of class $C^2$ and   if \eqref{eq:problemb} admits a solution.
System \eqref{eq:problemb}  was  considered in eulidean space by J. Serrin in  1971 in his seminal paper \cite{Serrin}.
%
\begin{Theorem}[Serrin \cite{Serrin}, Weinberger \cite{Wein}]\label{th:Serrin-Wein}
Bounded Serrin domains in  Euclidean space are  balls.
\end{Theorem}
  Serrin's argument   relies on  Alexandrov's  moving plane method,   with refined comparison principle. The proof of  Weinberger \cite{Wein} uses the  $P$-function method, Pohozaev identity  and maximum principles.
    Serrin's  result can be also derived from  Alexandrov's rigidity result. Namely if \eqref{eq:problemb}
  has a solution then $\de\Omega$ has CMC and thus is a ball, see  Farina and Kawohl in \cite{FK} and   Choulli and Henrot \cite{CH}. \\  
  The additional Neumann boundary condition in \eqref{eq:problemb} arises in many applications as a shape optimization problem for the underlying domain $\Omega$. For a detailed discussion of some applications e.g. in fluid dynamics and the linear theory of torsion,
see \cite{Serrin,Sirakov}.   Moreover,  as observed   in \cite{I.A.M}, Serrin domains arise also in the context of Cheeger sets in a Riemannian framework. To explain this connection more precisely, we recall that the Cheeger constant of a Lipschitz subdomain $\Omega \subset \cM$ is given by
$$
h(\O):= \inf_{A\subset \O }\frac{P(A)}{|A|}.
$$
Here the infimum is taken over measurable  subsets $A
\subset \Omega$, with finite perimeter $P(A)$ and where $|A|$  denotes the volume   of $A$  (both with respect to the metric $g$).  If $h(\Omega)$ is uniquely 
attained by $\O$ itself, then  $\Omega$ is   called   uniquely self-Cheeger. By means of the Weinberger's  $P$-function method, it is shown in \cite {I.A.M} that every bounded Serrin domain in a compact Riemannian manifold $\cM$, with   Ricci curvature bounded below by some constant, is uniquely self-Cheeger. Cheeger constants play an important role in eigenvalue estimates on Riemannian manifolds (see \cite{chavel}), whereas in the classical Euclidean case $(\cM,g)=(\R^N,g_{eucl})$ these notions have applications in the denoising  problem in  image processing, see e.g. \cite{EParini,Leon}.\\
The above Serrin's classification result was extended in \cite{Ku-Pra} to subdomains $\Omega$ of the round hemis-sphere $\cM = S^{N}$ or $\cM$ is a space forms of constant negative sectional curvatures. More precisely, it is proved in \cite{Ku-Pra} that any smooth Serrin domain $\Omega$ contained in a hemisphere of $S^N$ is a geodesic ball while the geodesic ball is the only Serrin domain in  space forms of constant negative sectional curvatures.
 It is an interesting and widely open problem to construct and classify Serrin domains. 
 
  Minlend and the author in \cite{FM},   found that on any compact Riemaniann manifold, there exists a Serrin  domain, which is a perturbation of a geodesic balls.  On the hand the symmetry group of the ambient manifold might be used to find nontrivial Serrin domains with different geometry.  Indeed, very recently, Minlend, Weth and the author in \cite{FMW-calc-var,FMW-ARMA} considered the cases $\cM=S^N$ the unit sphere and  the case $\M=\R^n\times \R^m/2\pi\Z^m$,  endowed with the flat metric, proving existences  of Serrin domains with nonconstant principal curvatures of the boundary. 

From an analytic point of view, the question of finding Serrin domains share similar features to  the one of finding CNMC hypersurfaces. Indeed,   an admissible class of parametrizations $\vp$ of the unknown domain boundary is considered. Then
 the solvability condition is formulated as an operator equation of the form $H(\vp)=const.$,
where $H$ is a nonlinear  Dirichlet-to-Neumann operator, thus sharing same traits with a  quasilinear nonlocal operators.    In fact, 
the nonlocal character of
this operator become apparent when    computing its linearization along a trivial branch of the problem, see e.g. \cite{FMW-calc-var,FMW-ARMA}.    Here we end up with a pseudo-differential operator of order 1 that  shares   many features with the square-root of the negative of the   Laplace operator, see \cite{GKS}.


From a geometric point of view, the results in \cite{FMW-calc-var,FMW-ARMA,CFW-2016, Cabre-Alex-Del-2015A} support naturally  the perception that the structure of the set of Serrin domains in a manifold $(\cM,g)$ share  similarites to those of CNMC hypersurfaces,  which  both,  up to a \textit{dimesion shift},  share  some   structure to   CMC hypersurfaces  in $\cM$.  Indeed, as in   the  theory  CNMC hypersurfaces, there exists in $\R^2$, non-straight periodic  Serrin domains, see \cite{FMW-ARMA}.  We also emphasize from \cite{FMW-calc-var} that there exist Serrin domains in $S^2$  which are not bounded by geodesic circles, whereas CMC-hypersurfaces in $S^2$ are obviously trivial, i.e., they are geodesic circles.\\
Next,   we recall that Alexandrov also proved in \cite{Alexandrov-62} that any closed embedded CMC hypersurface contained in a hemisphere of $S^N$ is a geodesic sphere. An explicit family of embedded  CMC hypersurfaces in $S^N$ with nonconstant principal curvatures was found by Perdomo in \cite{Perdomo} in the case $N \ge 3$. These hypersurfaces seem somewhat related to the Serrin domains  in \cite{FMW-calc-var}, although they bound a tubular neighborhood of $S^1$ and not of $S^{N-1}$. \\
 Overdetermined boundary value problems involving different elliptic equations has intensively studied recently. Starting with the pioneering papers of Pacard, Sicbaldi \cite{FPP},  Sicbaldi\cite{Sic} and Hauswirth, H\'elein and Pacard \cite{hauswirth-et-al}, the construction of nontrivial domains giving rise to solutions of overdetermined problems has been performed in many specific settings, see e.g. \cite{del-Pino-pacard-wei,ScSi,Ros-Ruiz-Sicbaldi-2016,morabito-sicbaldi}. Moreover, rigidity results for these domains were derived in \cite{ farina-valdinoci:2010-1,farina-valdinoci:2010-2, Ros-Sicbaldi,Ros-Ruiz-Sicbaldi-2015,Traizet}.


   \label{References}

\end{document}